\documentclass{amsart}
\usepackage{amssymb,latexsym}
\usepackage{bbm}
\usepackage{bm}
\usepackage{graphicx, pgf, pst-all}
\usepackage{mathrsfs}
\theoremstyle{plain}
\newtheorem{theorem}{Theorem}

\newtheorem{proposition}{Proposition}

\theoremstyle{definition}
\newtheorem{definition}{Definition}

\theoremstyle{remark}

\newtheorem{remark}{\bf Remark}
\numberwithin{equation}{section}
\def\Xint#1{\mathchoice
{\XXint\displaystyle\textstyle{#1}}%
{\XXint\textstyle\scriptstyle{#1}}%
{\XXint\scriptstyle\scriptscriptstyle{#1}}%
{\XXint\scriptscriptstyle\scriptscriptstyle{#1}}%
\!\int}
\def\XXint#1#2#3{{\setbox0=\hbox{$#1{#2#3}{\int}$ }
\vcenter{\hbox{$#2#3$ }}\kern-.6\wd0}}

\def\dashint{\Xint-}

\newcommand{\bd}{\operatorname{BD}}
\newcommand{\bv}{\operatorname{BV}}

\newcommand{\ld}{\operatorname{LD}}

\newcommand{\di}{\operatorname{div}}

\newcommand{\dif}{\operatorname{d}\!}

\newcommand{\R}{\mathbb{R}}

\newcommand{\locc}{\operatorname{loc}}

\newcommand{\dista}{\operatorname{dist}}

\newcommand{\sym}{\operatorname{sym}}

\newcommand{\trace}{\operatorname{Tr}}

\newcommand{\E}{\operatorname{E}\!}
\newcommand{\whr}{\operatorname{WRH}}

\newcommand{\ball}{\operatorname{B}}

\newcommand{\gm}{\operatorname{GM}}

\newcommand{\sobo}{\operatorname{W}}
\newcommand{\lebe}{\operatorname{L}}

\newcommand{\hold}{\operatorname{C}}

\newcommand{\D}{\operatorname{D}\!}

\renewcommand{\leq}{\leqslant}

\renewcommand{\di}{\operatorname{div}}
\newcommand{\bmo}{\operatorname{BMO}}
\newcommand{\sg}{\bm{\varepsilon}}

\begin{document}


\title[Regularity for Functionals of Linear Growth]
      {Symmetric--Convex Functionals of Linear Growth}
\author{Franz Gmeineder}
\address{Mathematical Institute\\
         University of Oxford\\
         Andrew Wiles Building\\
         OX2 6HG Oxford, United Kingdom} 
\email{gmeineder@maths.ox.ac.uk}
\urladdr{https://www.maths.ox.ac.uk/people/franz.gmeineder}
\date{\today}

\begin{abstract}
We discuss existence and regularity theorems for convex functionals of linear growth that depend on the symmetric rather than the full gradients. Due to the failure Korn's Inequality in the $\lebe^{1}$--setup, the full weak gradients of minima do not need to exist, and the paper aims for presenting methods that help to overcome these issues as to partial regularity and higher integrability of minimisers.
\end{abstract}

\maketitle

\section{Introduction}\label{intro} 
The purpose of the present paper is to survey and to announce existence and regularity results for minima of autonomous variational integrals which depend on the symmetric rather than the full gradients. More precisely, let $\Omega$ be an open and bounded Lipschitz subset of $\R^{n}$ and consider the variational principle
\begin{align}\label{eq:varprin}
\text{to minimise}\;\mathfrak{F}[v]:=\int_{\Omega}f(\sg(v))\dif x\;\text{over a Dirichlet class}\;\mathscr{D}, 
\end{align}
where $\sg(v):=\frac{1}{2}(\D v+\D^{\mathsf{T}}v)$ is the symmetric part of the weak gradient of a function $v\colon\Omega\to\R^{n}$, and $f\in\hold(\R_{\sym}^{n\times n})$ is a convex function of \emph{linear growth}. By the latter, we understand that there exist $c_{1},c_{2}>0$ such that 
\begin{align}\label{eq:lingrowth}
c_{1}|\mathbf{Z}|\leq f(\mathbf{Z})\leq c_{2}(1+|\mathbf{Z}|)\qquad\text{for all}\;\mathbf{Z}\in\R_{\sym}^{n\times n}. 
\end{align}
Under these conditions imposed on $f$, $\mathfrak{F}$ is well--defined on the space $\ld(\Omega)$ consisting of all $v\in\lebe^{1}(\Omega;\R^{n})$ whose weak symmetric gradients belong to $\lebe^{1}(\Omega;\R_{\sym}^{n\times n})$. This space is equipped with the canonical norm $\|v\|_{\ld}:=\|v\|_{\lebe^{1}}+\|\sg(v)\|_{\lebe^{1}}$, and we define $\ld_{0}(\Omega)$ to be the $\|\cdot\|_{\ld}$--closure of $\hold_{c}^{1}(\Omega;\R^{n})$. Consequently, we put $\mathscr{D}:=u_{0}+\ld_{0}(\Omega)$ for some fixed $u_{0}\in\ld(\Omega)$ and easily conclude that $\mathfrak{F}$ is coercive on $\mathscr{D}$ with respect to the $\ld$--norm. 

It is important to note that the aforementioned coerciveness fails when $\ld$ is replaced by $\sobo^{1,1}$. The reason for this is the lack of Korn's Inequality, a fundamental obstruction which we briefly describe now. Given $1<p<\infty$, Korn's Inequality asserts that there exists a constant $C>0$ such that 
\begin{align}\label{eq:Korn}
\int_{\Omega}|\D v|^{p}\dif x \leq C \int_{\Omega}|\sg(v)|^{p}\dif x
\end{align}
holds for all $v\in\hold_{c}^{1}(\Omega;\R^{n})$. In consequence, for convex and continuous integrands $g\colon \R_{\sym}^{n\times n}\to\R$ which satisfy $c_{1}|\mathbf{Z}|^{p}\leq g(\mathbf{Z})\leq c_{2}(1+|\mathbf{Z}|^{p})$ for all $\mathbf{Z}\in\R_{\sym}^{n\times n}$ and two constants $c_{1},c_{2}>0$, \eqref{eq:Korn} implies that the variational integral $\mathfrak{G}[v]:=\int_{\Omega}g(\sg(v))\dif x$ not only is well--defined on $\sobo^{1,p}(\Omega;\R^{n})$ but also coercive on suitable Dirichlet subclasses of $\sobo^{1,p}(\Omega;\R^{n})$ with respect to the usual $\sobo^{1,p}$--norm. As shall be explained in more detail in section \ref{sec:Ornstein} below, the reason for \eqref{eq:Korn} to hold is that the map $\Phi\colon \sg (v)\mapsto \D v$ (where $v\in\hold_{c}^{1}(\Omega;\R^{n})$ is tacitly identified with its trivial extension to the entire $\R^{n}$) is a singular integral of convolution type. Thus, by standard results for such operators, $\Phi$ is of strong--$(p,p)$ type if and only if $1<p<\infty$. In turn, if $p=1$, inequality \eqref{eq:Korn} fails to hold true, a fact which is often referred to as \emph{Ornstein's Non--Inequality}. Even stronger statements are available, some of which shall be discussed in section \ref{sec:Ornstein}. 

Since $\ld(\Omega)$ is a non--reflexive space, the second chief obstruction is that minimising sequences $(v_{k})\subset\mathscr{D}$ might not possess weakly convergent subsequences even though they are uniformly bounded with respect to the $\ld$--norm. To overcome this lack of compactness, it is reasonable to define the space $\bd(\Omega)$ of \emph{functions of bounded deformation} as the collection of all $v\in\lebe^{1}(\Omega;\R^{n})$ such that the distributional symmetric gradient $\sg(v)$ can be represented by a $\R_{\sym}^{n\times n}$--valued Radon measure of finite total variation on $\Omega$, in symbols $\sg(v)\in\mathcal{M}(\Omega;\R_{\sym}^{n\times n})$; see \cite{ST,ACD} for a detailled treatment of these spaces. In particular, by Ornstein's Non--Inequality, there exist elements $v\in\bd(\Omega)$ such that $\D v\notin\mathcal{M}(\Omega;\R^{n\times n})$ and hence $\bd(\Omega)$ contains $\bv(\Omega;\R^{n})$ as a proper subspace. In many respects, the properties of $\bd$--functions are reminiscent of those of $\bv$--functions, and we shall discuss similarities and discrepancies between the two function spaces as the paper evolves. Before passing to criteria that ensure the regularity of minima, we briefly revisit the treatment of the boundary value problem in $\bd$ which appears in a similar vein as that in $\bv$ as set up in the fundamental work of Giaquinta, Modica and Sou\v{c}ek \cite{GMS}.
\subsection{Relaxation and Generalised Minima}\label{sec:relaxed}
As an easy consequence of the Banach--Alaoglu and Rellich--Kondrachov Theorems, uniformly bounded sequences in $\ld(\Omega)$ possess subsequences that converge to some $v\in\bd(\Omega)$ in the weak*--sense. By this we understand that for some $(v_{j(k)})\subset (v_{j})$ there holds $v_{j(k)}\to v$ strongly in $\lebe^{1}(\Omega;\R^{n})$ and $\sg(v_{j(k)})\stackrel{*}{\rightharpoonup}\sg(v)$ in the sense  of weak*--convergence of Radon measures on $\Omega$ as $k\to\infty$. In this situation, the weak*--limit map can be shown to exist, however, to establish a reasonable notion of minimality for $v$, the functional $\mathfrak{F}$ must be extended to $\bd(\Omega)$ first. To keep the presentation simple, we stick to the $\lebe^{1}$--Lebesgue--Serrin extension given by 
\begin{align*}
\overline{\mathfrak{F}}[v]:=\inf\big\{\liminf_{k\to\infty}\mathfrak{F}[v_{k}]\colon\; (v_{k})\subset\mathscr{D},\;v_{k}\to v\;\text{in}\;\lebe^{1}(\Omega;\R^{n})\;\text{as}\;k\to\infty \big\},\qquad v\in\bd(\Omega). 
\end{align*}
Note that this type of relaxation is reasonable indeed: If $(v_{k})\subset\mathscr{D}$ is bounded with respect to the $\ld$--norm and converges to some $v\in\lebe^{1}(\Omega;\R^{n})$ strongly in $\lebe^{1}(\Omega;\R^{n})$, then $v\in\bd(\Omega)$ by lower semicontinuity of the total deformation $|\sg(\cdot)|(\Omega)$ with respect to strong $\lebe^{1}$--convergence. It needs to be noted that the functionals $\overline{\mathfrak{F}}$ admit the explicit integral representation 
\begin{align*}
\overline{\mathfrak{F}}[v]=\int_{\Omega}f(\mathscr{E}v)\dif x + \int_{\Omega}f^{\infty}\left(\frac{\dif\E v}{\dif |\E^{s}v|}\right)\dif |\E^{s}u| + \int_{\partial\Omega}f^{\infty}(\trace(v-u_{0})\odot\nu_{\partial\Omega})\dif\mathcal{H}^{n-1}
\end{align*}
for $v\in\bd(\Omega)$, where 
\begin{align*}
\sg(v)= \E^{ac}v + \E^{s}v = \frac{\dif\E v}{\dif\mathscr{L}^{n}}\mathscr{L}^{n}+\frac{\dif \E^{s}v}{\dif |\E^{s}v|}|\E^{s}v| = \mathscr{E}v\mathscr{L}^{n}+\frac{\dif \E^{s}v}{\dif |\E^{s}v|}|\E^{s}v|
\end{align*}
is the Radon--Nikod\v{y}m decomposition of $\sg(v)$ into its absolutely continuous and singular parts with respect to Lebesgue measure $\mathscr{L}^{n}$; moreover, $f^{\infty}\colon\R_{\sym}^{n\times n}\to\R$ is the \emph{recession function} of $f$ defined as
\begin{align*}
f^{\infty}(\mathbf{Z}):=\lim_{t\searrow 0}tf\left(\mathbf{Z}/t\right),\qquad\mathbf{Z}\in\R_{\sym}^{n\times n}, 
\end{align*}
and captures the behaviour of the integrand at infinity, that is to say, where the Lebesgue density of $\sg(v)$ with respect to $\mathscr{L}^{n}$ becomes singular. Note that by \cite{ACD,Hajlasz}, the density $\tfrac{\dif\E v}{\dif\mathscr{L}^{n}}$ can be shown to equal the symmetric part $\mathscr{E}v$ of the approximate gradient of $v$ $\mathscr{L}^{n}$--a.e.. The trace terms in fact are sensible, as by \cite{ST,Baba}, $\bd$--functions attain boundary values in the $\lebe^{1}$--sense. Noting that for $a,b\in\R^{n}$, $a\odot b:=\tfrac{1}{2}(ab^{\mathsf{T}}+ba^{\mathsf{T}})$ is the symmetric tensor product and $\nu_{\partial\Omega}$ the outer unit normal to $\partial\Omega$, the boundary integral term appearing in the integral representation admits the interpretation of a penalisation term that leads to larger values of the functional provided the $\lebe^{1}$--distance of $\trace(v)$ from $\trace(u_{0})$ is increased. The proof of the above representation follows along the lines of \cite{Bild1} in the full gradient case, and for more background information, the reader is referred to \cite{GK}. For completeness, we make the following 
\begin{definition}[$\bd$--Minima]
An element $u\in\bd(\Omega)$ is called a \emph{$\bd$--minimiser} if and only if $\overline{\mathfrak{F}}[u]\leq\overline{\mathfrak{F}}[v]$ for all $v\in\bd(\Omega)$. 
\end{definition}
On a sidenote, let us remark that convexity of $f$ substantially simplifies the proof of the integral representation for the relaxed functional. In fact, in the convex setup, it is possible to use the Goffman--Serrin relaxation machinery \cite{GoffSerrin}, whereas in the quasiconvex situation more subtle arguments need to be invoked; see the work of Rindler \cite{Rindler} for more detail. 

Another notion of minimisers has been employed by Bildhauer \& Fuchs \cite{Bild2,Bild1} in the setting of linear growth functionals on $\bv$, whose adaption to the present situation reads as follows: 
\begin{definition}[Generalised Minima]\label{def:genmin}
Let $\Omega$ be an open and bounded Lipschitz subset of $\R^{n}$ and fix a boundary datum $u_{0}\in\ld(\Omega)$. The set of \emph{generalised minima} of $\mathfrak{F}$ given by \eqref{eq:varprin} consists of all those $u\in\bd(\Omega)$ for which there exists an $\mathfrak{F}$--minimising sequence $(u_{k})\subset \mathscr{D}_{u_{0}}:=u_{0}+\ld_{0}(\Omega)$ such that $u_{k}\to u$ strongly in $\lebe^{1}(\Omega;\R^{n})$ as $k\to\infty$. The set of all generalised minima is denoted $\gm(\mathfrak{F})$.
\end{definition}
Now, if $f\in\hold(\R_{\sym}^{n\times n})$ is convex -- that is to say, $f$ is \emph{symmetric--convex} -- then $u\in\bd(\Omega)$ can be shown to be a $\bd$--minimiser if and only if it is a generalised minimiser, and in this case there holds
\begin{align}
\mathfrak{F}[u]=\min \overline{\mathfrak{F}}(\bd(\Omega))=\inf\mathfrak{F}[\mathscr{D}]. 
\end{align}
The crucial point thus is to establish existence of a $\bd$--minimiser. This is, however, easily achieved by employing the direct method and by use of Reshetnyak--type theorems on the lower semicontinuity of functionals of measures with respect to the weak*-- and strict topologies; see \cite{Reshetnyak,Bild1}. Hence a satisfactory existence theory is established, and the foremost aim of the present paper is to survey the reguarity properties of generalised minima. 
\subsection{Organisation of the Paper and Description of Results}
Having settled existence of generalised minima in the previous section, the paper focusses on regularity results in all of what follows. In section \ref{sec:Ornstein}, we revisit Korn's Inequality and Ornstein's Non--Inequality and strengthen the sketchy arguments outlined above to understand the chief obstructions for regularity results in the symmetric--convex case. In section \ref{sec:regularity}, we report on recent developments regarding the H\"{o}lder and Sobolev regularity of generalised minima of symmetric--convex functionals subject to strong convexity conditions imposed on the variational integrands $f$. Firstly describing a partial regularity result due to the author \cite{Gm} in the spirit of Anzellotti \& Giaquinta \cite{AG}, we then turn to conditions on the variational integrands to produce generalised minima of class $\bv_{\locc}$ or even $\sobo_{\locc}^{1,p}$ for some $1<p<\infty$, the latter being joint work with Jan Kristensen \cite{GK}. To our best knowledge, these results are the first of their kind and extend the regularity theory on the Dirichlet problem on $\bv$ to that on $\bd$; see \cite{Bild1,BS1}. To conclude with, in section \ref{sec:Aconvex} we introduce the spaces $\sobo^{\mathbb{A},1}$ and $\bv_{\mathbb{A}}$ as suitable generalisatons of $\bv$ and $\bd$ and highlight open questions that would lead to a satisfactory existence and regularity theory in this fairly general setup. 
\section*{Acknowledgment}
The author is indebted to the University of Z\"{u}rich for financial support to attend the 9th European Conference on Elliptic and Parabolic Problems held in Gaeta in May 2016. Moreover, he gratefully acknowledges the comments of an anonymous referee which helped to improve the exposition of the material.
\section*{Notation}
The symmetric $n\times n$--matrices are denoted $\R_{\sym}^{n\times n}$, and we use the symbol $\langle\cdot,\cdot\rangle$ for the euclidean inner product on finite dimensional spaces. Lastly, $\mathscr{L}^{n}$ and $\mathcal{H}^{n-1}$ denote the $n$--dimensional Lebesgue-- or $(n-1)$--dimensional Hausdorff measures, respectively, and we use $(u)_{U}:=\dashint_{U}u\dif x:=(\mathscr{L}^{n}(U))^{-1}\int_{U}u\dif x$ for the mean value of a locally integrable function $u\colon U\to\R^{N}$ whenever this is well--defined.
\section{Korn's Inequality and Ornstein's Non--Inequality}\label{sec:Ornstein}
Before we embark on the regularity of generalised minima as addressed in the introduction, we briefly wish to comment on Korn's Inequality in slightly more detail. To this end, let $\mathbb{A}[D]$ be a linear, homogeneous first order and constant coefficient differential operator between the two finite--dimensional real vector spaces $V$ and $W$, i.e., $\mathbb{A}[D]$ can be written in the form 
\begin{align}\label{eq:differentialoperator}
\mathbb{A}[D]=\sum_{|\alpha|=1}\mathbb{A}_{\alpha}\partial^{\alpha}, 
\end{align}
where $\mathbb{A}_{\alpha}\colon V\to W$ are fixed linear mappings. We associate with $\mathbb{A}[D]$ its \emph{symbol map}
\begin{align}
\mathbb{A}[\xi]:=\sum_{|\alpha|=1}\xi_{\alpha}\mathbb{A}_{\alpha},\qquad\xi=(\xi_{1},...,\xi_{n})\in\R^{n}, 
\end{align}
and call $\mathbb{A}[D]$ \emph{elliptic} if and only if $\mathbb{A}[\xi]\colon V\to W$ is injective for any $\xi\neq 0$. Given an elliptic operator $\mathbb{A}[D]$ and $u\in\hold_{c}^{\infty}(\R^{n};\R^{n})$, we can thus retrieve $u$ from $\mathbb{A}[D]u$ by means of the operator
\begin{align*}
u(x)=\mathbf{G}[\mathbb{A}[D]u](x)=c_{n}\mathscr{F}_{\xi\mapsto x}^{-1}((\mathbb{A}^{*}[\xi]\circ\mathbb{A}[\xi])^{-1}\mathbb{A}^{*}[\xi]\widehat{\mathbb{A}[D]u})=:\Phi(\mathbb{A}[D]u)(x),\qquad x\in\R^{n}, 
\end{align*}
where $c_{n}>0$ is a constant and $\mathbb{A}^{*}[\xi]$ is the adjoint symbol of $\mathbb{A}[\xi]$ being defined in the obvious manner. Since $(\mathbb{A}^{*}[\xi]\circ\mathbb{A}[\xi])^{-1}\mathbb{A}^{*}[\xi]$ is homogeneous of degree $-1$, it is easily seen that $\mathbf{G}$ is a Riesz potential operator of order $1$, in particular, we have the bound 
\begin{align*}
|\mathbf{G}[v](x)|\lesssim \int_{\R^{n}}\frac{|v(y)|}{|x-y|^{n-1}}\dif y\qquad\text{for all}\;x\in\R^{n}.
\end{align*}
Finally, differentiating $\mathbf{G}[\mathbb{A}[D]u]$ immediately yields that $\D u$ can be written as a singular integral of convolution type acting on $\mathbb{A}[D]u$; see \cite{Stein}. Therefore, the operator $\Phi\colon \mathbb{A}[D]u\mapsto \D u$ given by $\Phi (\mathbb{A}[D]u):=\D\big(\mathbf{G}[\mathbb{A}[D]u]\big)$ extends to a bounded linear operator $\widetilde{\Phi}\colon\lebe^{p}(\R^{n};W)\to\lebe^{p}(\R^{n};\R^{n}\times V)$ provided $1<p<\infty$. Hence, given $1<p<\infty$ and an open subset $\Omega$ of $\R^{n}$, Korn's inequality $\|\D u\|_{\lebe^{p}(\Omega;\R^{n}\times V)}\leq C\|\mathbb{A}[D]u\|_{\lebe^{p}(\Omega;W)}$ for all $u\in\hold_{c}^{\infty}(\Omega;V)$ with a finite constant $C=C(\mathbb{A},p)>0$ follows from the aforementioned boundedness of singular integrals by extending elements of $\hold_{c}^{\infty}(\Omega;V)$ to $\R^{n}$ by zero. 

Korn's Inequality in the form as given above can be generalised to various other settings; see \cite{MNRR} and \cite{BD,BCD} for more recent developments in the context of Orlicz functions. Its failure in the case $p=1$ has witnessed a variety of notable contributions beyond Ornstein's original article \cite{Ornstein}; see, among others, \cite{CFMM,KirchKrist}. In particular, as can be seen best through Theorem 1.3 of Kirchheim and Kristensen's study \cite{KirchKrist}, there are only trivial $\lebe^{1}$--estimates in the following sense:
\begin{theorem}[\cite{KirchKrist}, Theorem 1.3]\label{thm:Ornstein}
Let $V,W,X$ be three finite--dimensional vector spaces and consider two $k$--th order linear and homogeneous differential operators $\mathbb{A}_{1}[D]$ and $\mathbb{A}_{2}[D]$ of the form 
\begin{align*}
\mathbb{A}_{1}[D]=\sum_{|\alpha|=k}\mathbb{A}_{\alpha}^{1}(x)\partial^{\alpha}\quad\text{and}\quad\mathbb{A}_{2}[D]=\sum_{|\alpha|=k}\mathbb{A}_{\alpha}^{2}(x)\partial^{\alpha}, 
\end{align*}
with locally integrable coefficients $\mathbb{A}_{\alpha}^{1}\in\lebe_{\locc}^{1}(\R^{n};\mathscr{L}(V;W))$ and $\mathbb{A}_{\alpha}^{2}\in\lebe_{\locc}^{1}(\R^{n};\mathscr{L}(V;X))$ for all $|\alpha|=k$, respectively, the following are equivalent: 
\begin{enumerate}
\item There exists a constant $C>0$ such that $\|\mathbb{A}_{2}[D]\varphi\|_{\lebe^{1}(\R^{n};W)}\leq c\|\mathbb{A}_{1}[D]\varphi\|_{\lebe^{1}(\R^{n};W)}$ holds for all $\varphi\in\hold_{c}^{\infty}(\R^{n};V)$. 
\item There exists $C\in\lebe^{\infty}(\R^{n};\mathscr{L}(W;X))$ with $\|C\|_{\lebe^{\infty}(\R^{n};\mathscr{L}(W;X))}\leq c$ such that $\mathbb{A}_{\alpha}^{2}(x)=C(x)\mathbb{A}_{\alpha}^{1}(x)$ for $\mathscr{L}^{n}$--a.e. $x\in\R^{n}$ and each $\alpha\in\mathbb{N}_{0}^{n}$ with $|\alpha|=k$. 
\end{enumerate}
\end{theorem}
Here, triviality of $\lebe^{1}$--estimates means that \emph{if} a Korn--type inequality holds (1), then the coefficients are the same up to multiplication with an $\lebe^{\infty}$--function. It is important to note that both the symmetric gradient operator $\mathbb{A}_{1}[D]u=\sg(u)$ or the trace--free symmetric gradient operator $\mathbb{A}_{1}[D]u=\sg(u)-\frac{1}{n}\di(u) \mathbbm{1}_{n\times n}$ with the $(n\times n)$--unit matrix $\mathbbm{1}_{n\times n}\in\R^{n\times n}$ do \emph{not} verify (2) with $\mathbb{A}_{2}[D]=\D u$ and hence they do not admit Korn--type inequalities in the $\lebe^{1}$--setup. 
\section{Partial and Sobolev Regularity}\label{sec:regularity}
Besides existence of generalised minima as outlined in the introduction, it is natural to investigate their regularity properties, a program which has been launched in the $\bv$--setting in \cite{GMS,LU}. Here, we focus on H\"{o}lder-- and Sobolev regularity and shall describe the main obstructions that come along both with the linear growth hypothesis and Ornstein's Non--Inequality. In particular, the latter motivates to study conditions imposed on the variational integrand $f$ that guarantee existence of the full gradients of generalised minima as elements of $\mathcal{M}(\Omega;\R^{n\times n})$ or $\lebe^{1}(\Omega;\R^{n\times n})$. 
\subsection{H\"{o}lder Regularity}
To begin with, let us note that by the genuine vectorial nature of the functional $\mathfrak{F}$, $\bd$--minima cannot be shown to share everywhere $\hold^{1,\alpha}$--regularity unless strong structural conditions are imposed on the variational integrands $f$. This is in the spirit of the famous counterexamples due to De Giorgi \cite{DG} and Giusti \& Miranda \cite{GiuMir} (see also \cite{Mingione} for an excellent overview) which demonstrate that in the case $N>1$, functionals of the form 
\begin{align*}
\mathfrak{G}[v]:=\int_{\Omega}g(\nabla v)\dif x,\qquad v\colon \Omega\to\R^{N}
\end{align*}
do not necessarily produce minimisers of class $\hold_{\locc}^{1,\alpha}(\Omega;\R^{N})$ even if suitable ellipticity, boundedness and measurability assumptions are imposed on the integrands $g$. The correct substitute is then given by the notion of \emph{partial regularity}: Given $u\in\bd(\Omega)$, we define its \emph{regular set}
\begin{align}
\Omega_{u}:=\{x\in\Omega\colon\; \nabla u\;\text{is H\"{o}lder continuous in a neighbourhood of}\;x\}, 
\end{align}
and note that the definition of $\Omega_{u}$ depends on the full weak gradients $\nabla u$ indeed. We call $\Omega_{u}$ the \emph{regular set} of $u$ and its relative complement $\Sigma_{u}:=\Omega\setminus\Omega_{u}$ the \emph{singular set}. Adopting these notions, we say that $u\in\bd(\Omega)$ is \emph{partially regular} if and only if $\Omega_{u}$ is open and $\mathscr{L}^{n}(\Sigma_{u})=0$. 

Referring the reader to \cite{Beck,Giusti,Mingione} for a comprehensive overview of techniques to establish partial regularity of minima of elliptic variational integrals, we note that most approaches to the partial regularity rely on the higher integrability of gradients. Such higher integrability results in turn often stem from Caccioppoli--type inequalities in conjunction with Gehring's lemma. Indeed, to sketch the prototype form of such an argument, let $1<p<2$ and assume that $u\in\sobo^{1,p}(\Omega;\R^{N})$ satisfies a Caccioppoli--type inequality of the form
\begin{align}\label{eq:Cacc}
\dashint_{\ball(z,r)}|\D u|^{p}\dif x \leq C\dashint_{\ball(z,2r)}\left\vert \frac{u-(u)_{\ball(z,2r)}}{r}\right\vert^{p}\dif x
\end{align}
for all $z\in\Omega$ and $0<r<\dista(z,\partial\Omega)/2$. Now, applying the Sobolev--Poincar\'{e}--inequality to the right side, we deduce that $\dashint_{\ball(z,r)}|\D u|^{p}$ can be locally estimated against $\dashint_{\ball(z,2r)}|\D u|^{q}\dif x$ for some $1<q<p$, and in this sense $\D u$ satisfies a reverse H\"{o}lder inequality. By Gehring's Lemma, we then obtain that $\D u$ belongs to some $\lebe_{\locc}^{p+\varepsilon}$ for some $\varepsilon>0$, and the reader will notice that the above argument remains unchanged when $\D$ is replaced by $\sg$ and the mean values on the right side of \eqref{eq:Cacc} by suitable rigid deformations, that is, elements of the nullspace of $\sg$. 

If $p=1$, then Gehring's Lemma self--improves the integrability of $\D u$ or $\sg(u)$, respectively, if and only if we can choose $q<p=1$. Thus a suitable Sobolev--Poincar\'{e} inequality would be required, estimating the $\lebe^{1}$--norm of a function against the $\lebe^{q}$--norm of its gradient or symmetric gradient, respectively. However, in the full generality as needed here, this is ruled out by a counterexample due to Buckley \& Koskela \cite{BuckleyKoskela}. In turn, one is lead to the so--called \emph{weak reverse H\"{o}lder classes} $\whr$ whose elements satisfy suitable Sobolev--Poincar\'{e} inequalities for $q<p=1$ by definition in a natural way; see \cite{BuckleyKoskela} for more background information. Coming back to the higher integrability addressed above, it is not clear that the gradients of minima belong to such weak reverse H\"{o}lder classes at the relevant stage of the proof; hence different methods are required in the linear growth setting. In the symmetric--convex case as described in the introduction, such estimates can be achieved by use difference quotient--type methods, but in turn require strong ellipticity assumptions on the integrands; see section \ref{sec:Sobolev} below. 

A direct approach to the partial regularity that is particularly designed for convex functionals and applies to functionals of linear growth, too, is that of Anzellotti \& Giaquinta \cite{AG}. As usual, this particular method also relies on decay estimates for suitable excess functionals too. To describe the decisive feature of this method, let us remark that unlike other, perhaps more standard approaches like blow--up proofs, Anzellotti \& Giaquinta derive the required decay estimates through comparing minima with suitable mollifications thereof. These mollifications in turn are shown to be close to solutions of elliptic second order PDE and thus enjoy good decay estimates which then are shown to carry over to the generalised minima themselves. Relying on mollifications and, consequently, Jensen's inequality, the method is well--suited for convex problems, whereas it is not clear how to generalise it to quasiconvex integrands, for instance. With the case of full gradients being treated in \cite{AG}, the respective generalisation to the symmetric gradient case will be given in \cite{Gm} by the following
\begin{theorem}
Let $f\in \hold^{2}(\R_{\sym}^{n\times n};\R_{\geq 0})$ be convex and of linear growth. Suppose that $u\in\bd(\Omega)$ is a $\bd$--minimiser of $\mathfrak{F}$ given by \eqref{eq:varprin}. If $(x,z)\in\Omega\times\R_{\sym}^{n\times n}$ is such that 
\begin{align*}
\lim_{R\searrow 0}\left[\dashint_{\ball(x,R)}|\mathscr{E}u-z|\dif x + \frac{|\E^{s}u|(\ball(x,R))}{\mathscr{L}^{n}(\ball(x,R))}\right]=0
\end{align*}
and $f''(z)$ is positive definite, then $u\in \hold^{1,\alpha}(U;\R^{n})$ for a suitable neighbourhood $U$ of $x$ for all $0<\alpha<1$. 
\end{theorem}
Assuming the theorem, the standard Lebesgue differentiation theorem for Radon measures yields the claimed partial regularity: There exists an open subset $\Omega_{u}$ of $\Omega$ with $\mathscr{L}^{n}(\Omega\setminus\Omega_{u})=0$ such that for any $x\in\Omega_{u}$ there exists $r>0$ with $u\in\hold^{1,\alpha}(\ball(x,r);\R^{n})$ for every $0<\alpha<1$. We wish to conclude with the following
\begin{remark}[Non--Autonomous Integrands]
Since $\mathfrak{F}$ given by \eqref{eq:varprin} is autonomous, it is possible to overcome the higher integrability of the symmetric gradients in the proof of the above theorem. If the integrand in addition is $x$--dependent, then the higher integrability seems to be a necessary to conclude the result in this non--autonomous case too. Going back to the discussion at the beginning of the section, such a result is therefore unlikely to be established by means of the method as described above.

\end{remark}
\subsection{Sobolev Regularity}\label{sec:Sobolev}
Next we turn to Sobolev regularity of generalised minima and hereafter aim for conditions on the integrands $f$ under which generalised minima genuinely belong to $\bv_{\locc}(\Omega;\R^{n})$ or $\sobo_{\locc}^{1,p}(\Omega;\R^{n})$ for some $1<p<\infty$. To obtain such results, we shall work with a strong convexity adapted from that of Bildhauer \& Fuchs \cite{Bild2} in the full gradient case:
\begin{definition}[$\mu$--ellipticity]\label{def:muelliptic}
Let $\mu>1$. A $\hold^{2}$--integrand $f\colon \R_{\sym}^{n\times n}\to \R_{\geq 0}$ is called $\mu$--\emph{elliptic} if and only if there exist $0<\lambda\leq\Lambda<\infty$ such that 
\begin{align}\label{eq:muell}
\lambda\frac{|\mathbf{A}|^{2}}{(1+|\mathbf{B}|^{2})^{\frac{\mu}{2}}}\leq \langle f''(\mathbf{B})\mathbf{A},\mathbf{A}\rangle \leq \Lambda \frac{|\mathbf{A}|^{2}}{(1+|\mathbf{B}|^{2})^{\frac{1}{2}}}
\end{align}
holds for all $\mathbf{A},\mathbf{B}\in \R_{\sym}^{n\times n}$. We further say that the variational integral $\mathfrak{F}$ is $\mu$--elliptic provided its integrand $f$ is. 
\end{definition} 
Prototypical examples are given by the area--type integrands $m_{p}(\bm{\xi}):=(1+|\bm{\xi}|^{p})^{\frac{1}{p}}$ for $\bm{\xi}\in\R_{\sym}^{n\times n}$, $1<p<\infty$; so, for instance, the area integrand $m_{2}$ is $3$--elliptic, whereas $m_{p}$ coincides with a $\mu=p+1$--elliptic integrand away from the unit ball; also see \cite{Bild2} for more detail. 
It is important to remark that $\mu=1$ is explicitely excluded in definition \ref{def:muelliptic}; indeed, $1$--elliptic integrands correspond to $L\log L$--growth. For such integrands, $\sg(u)\in L\log L_{\locc}(\Omega;\R_{\sym}^{n\times n})$ already implies $\D u\in\lebe_{\locc}^{1}(\Omega;\R^{n\times n})$ and hence the full gradients are known to exist and belong to $\lebe^{1}$ locally. In this setup of $L\log L$--growth, the corresponding regularity theory has been established by Fuchs, Seregin and collaborators \cite{FrehseSeregin,FuchsSeregin1,FuchsMingione} among others; see the extensive monograph \cite{FuchsSeregin2} for more information. 
\subsubsection{Results on the Dirichlet Problem on $\bv$}
To explain our method, it is convenient to firstly report on the available higher integrability results for $\mu$--elliptic functionals in the full gradient case 
\begin{align}
\mathcal{F}[u]:=\int_{\Omega}f(\nabla v)\dif x,\qquad v\in\mathcal{D}:=u_{0}+\sobo_{0}^{1,1}(\Omega;\R^{N}), 
\end{align}
where $u_{0}\in\sobo^{1,1}(\Omega;\R^{N})$ is a given boundary datum. In analogy with Definition \ref{def:genmin}, we say that $v\in\bv(\Omega;\R^{N})$ is a \emph{generalised minimiser} for $\mathcal{F}$ if and only if there exists an $\mathcal{F}$--minimising sequence $(v_{k})\subset\mathcal{D}$ such that $v_{k}\to v$ strongly in $\lebe^{1}(\Omega;\R^{N})$ as $k\to\infty$. 

Using a vanishing viscosity approach, Bildhauer \cite{Bild3} provided the first higher integrability results for gradients of generalised minima. Precisely, assuming $u_{0}\in\sobo^{1,2}(\Omega;\R^{N})$ for the boundary data, the functional $\mathcal{F}$ is stabilised by adding small Laplacians, i.e., we consider
\begin{align*}
\mathcal{F}_{\delta}[v]:=\mathcal{F}[v]+\frac{\delta}{2}\int_{\Omega}|\nabla v|^{2}\dif x\qquad\text{on}\;\;\mathcal{D}:=u_{0}+\sobo_{0}^{1,2}(\Omega;\R^{N})
\end{align*}
and finally aim for sending $\delta\searrow 0$. Denoting the unique minimiser of $\mathfrak{F}$ over $\mathcal{D}$ by $u_{\delta}$, it is easy to prove that $(u_{\delta})$ is a minimising sequence for $\mathcal{F}$. Bildhauer, in turn building on ideas of Seregin \cite{Seregin}, then was able to show that if $(u_{\delta})$ satisfies the \emph{local boundedness assumption}
\begin{align}\label{eq:locbound}
\text{for all}\;K\Subset \Omega\;\text{there exists}\;C(K)>0\;\text{with}\;\sup_{0<\delta<1}\|u_{\delta}\|_{\lebe^{\infty}(K;\R^{N})}\leq C(K), 
\end{align}
then the weak*--limit $u$ of $(u_{\delta})$ belongs to $\sobo_{\locc}^{1,p}(\Omega;\R^{N})$ for some $p=p(\mu)>1$ provided $1<\mu<3$, and to $\sobo_{\locc}^{1,L\log L}(\Omega;\R^{N})$ provided $\mu=3$. Apart from the strong assumptions made on the particular minimising sequence, the boundary data and the functional itself, the result merely applies to one particular such generalised minimiser. The reason for this is the possible non--uniqueness of generalised minima which, in turn, is a consequence of the recession parts in the relaxed variational integral. Indeed, even if $f\in\hold^{2}(\R^{N\times n})$ is strictly convex, the recession function $f^{\infty}\colon\R^{N\times n}\to\R$ is positively $1$--homogeneous and thus never strictly convex. In consequence, if a minimiser does not have vanishing singular part with respect to Lebesgue measure, uniqueness in general fails as $f$ and $f^{\infty}$ act on two mutually singular parts of the gradients. To achieve uniqueness, one must therefore genuinely rule out the singular parts of minima. This has been achieved recently by Beck \& Schmidt \cite{BS1} by sophisticated use of the Ekeland variational principle in the negative Sobolev space $\sobo^{-1,1}$ (see Prop. \ref{prop:ekeland} below) for the borderline case $\mu=3$. Referring the reader for the precise outline to \cite{BS1}, the general streamline is this: Starting from an arbitrary minimising sequence $(u_{k})\subset u_{0}+\sobo_{0}^{1,1}(\Omega;\R^{N})$, the Ekeland variational principle yields another minimising sequence $(v_{k})$ that is $\sobo^{-1,1}$--close to $(u_{k})$ and has the same weak*--limit. This new sequence $(v_{k})$ is then shown to be a sequence of almost minimisers to suitably stabilised functionals and thus can be proved to belong to $\sobo_{\locc}^{2,2}$. At this point it is possible to adapt Bildhauer's approach and hence, by arbitrariness of $(u_{k})$, uniqueness and the aforementioned higher regularity results follow at once. However, it needs to be stressed that Beck \& Schmidt's method of proof also relies on a version of the local boundedness assumption. Further, assuming Uhlenbeck, i.e., radial structure of the integrands, stronger results such as $\hold^{1,\alpha}$--regularity of generalised minima of the $\mathcal{F}$ can be achieved; see \cite{Bild1,Bild3,BS2}. 
\subsubsection{Results on the Dirichlet Problem on $\bd$}
Going back to functionals $\mathfrak{F}$ as given by \eqref{eq:varprin}, the main difficulty lies in the appearance of the \emph{full difference quotients} when aiming for higher integrability estimates and testing the Euler--Lagrange equation of a suitably stabilised functional with the canonical choice $\varphi:=\Delta_{s,h}^{-}(\rho^{2}\Delta_{s,h}^{+}v)$, where $\Delta_{s,h}^{\pm}v(x):=\frac{1}{h}(v(x\pm he_{s})-v(x))$ are the forward or backward difference quotients, respectively. By Ornstein's Non--Inequality, $\Delta_{s,h}^{+}v$ cannot even be bounded locally in $\lebe^{1}$ for $v\in\bd$ in general, and thus the suitable device hence is to work with finite differences instead of difference quotients and to establish estimates for carefully chosen Besov--norms of the symmetric gradients. 

We pass on to a more precise description of the method which is, to some extent, inspired by \cite{BS1}. Let $(v_{k})\subset\mathscr{D}$ be a minimising sequence for the $\mu$--elliptic functional $\mathfrak{F}$ given by \eqref{eq:varprin}, with $\mu$ to be determined later on. Then we consider for a suitable sequence $(\alpha_{k})\subset\R_{>0}$ with $\alpha_{k}\searrow 0$ as $k\to\infty$ the stabilised functionals 
\begin{align*}
\mathfrak{F}_{k}[w]:=\int_{\Omega}f(\sg(w))\dif x+\alpha_{k}\int_{\Omega}(1+|\sg(w)|^{p})\dif x =:\int_{\Omega}f_{k}(\sg(w))\dif x
\end{align*}
with $p\geq n$ on appropriately modified Dirichlet classes $\mathscr{D}_{k}$ to keep track of the fact that the leading part of $\mathfrak{F}_{k}$ is the $p$--th Dirichlet energy; indeed, as $p>1$, minima of $\mathfrak{F}_{k}$ belong to $\sobo^{1,p}$ and thus possess full gradients in $\lebe^{p}$ by Korn's Inequality. Extending each $\mathfrak{F}_{k}$ to $(\sobo_{0}^{1,\infty}(\Omega;\R^{n}))^{*}$ by infinity on $(\sobo_{0}^{1,\infty}(\Omega;\R^{n}))^{*}\setminus\mathscr{D}_{k}$, we obtain a lower semicontinuous functional on $(\sobo_{0}^{1,\infty}(\Omega;\R^{n}))^{*}$ which is continuous with respect to the norm topology on $(\sobo^{1,\infty}(\Omega;\R^{n}))^{*}$. To continue, we recall the following instrumental
\begin{proposition}[Ekeland's Variational Principle]\label{prop:ekeland}
Let $(X,d)$ be a complete metric space and assume that $f\colon X\to [0,\infty]$ is  lower semicontinuous with $\inf F(X)<\infty$. If for some $\varepsilon>0$ and $x\in X$ there holds $F[u]\leq \inf F(X)+\varepsilon$, then there exists $v\in X$ such that $d(x,v)\leq\sqrt{\varepsilon}$ and 
\begin{align*}
F[v]\leq F[w]+\sqrt{\varepsilon}d(v,w)\qquad\text{for all}\;w\in X. 
\end{align*}
\end{proposition}
For a proof and a discussion of this result, see \cite{Giusti}, Thm. 5.6.. Using suitable approximations and Proposition \ref{prop:ekeland} in the Banach space $(\sobo_{0}^{1,\infty})^{*}$, we obtain a sequence $(u_{k})$, each of whose members is an almost minimiser of $\mathfrak{F}$, is $(\sobo_{0}^{1,\infty})^{*}$--close to $v_{k}$ and, most crucially, each $u_{k}$ satisfies the perturbed Euler--Lagrange equation 
\begin{align*}
\left\vert\int_{\Omega}\langle f'_{k}(\sg(u_{k})),\sg(\varphi)\rangle\dif x\right\vert \leq \frac{1}{k}\|\varphi\|_{(\sobo_{0}^{1,\infty}(\Omega;\R^{n})^{*}}+\big(\text{small perturbation}\big)
\end{align*}
for all $\varphi\in\sobo_{0}^{1,p}(\Omega;\R^{n})$. Given $x_{0}\in\Omega$, $0<r<R<\dista(x_{0},\partial\Omega)$, we consider for an arbitrary standard unit vector $e_{s}$, $s=1,...,n$, the test functions $\varphi:=\tau_{s,h}^{-}(\rho^{2}\tau_{s,h}^{+}u_{k})$ with $\tau_{s,h}^{\pm}=h\Delta_{s,h}^{\pm}$. Essentially following, e.g., \cite{Giusti}, section 10.1,  and employing $\mu$--ellipticity of $f$, we end up with a coercive inequality
\begin{align}\label{eq:coercive}
\int_{\Omega} \frac{|\rho\tau_{s,h}^{+}\sg(u_{k})|^{2}}{(1+|\sg(u_{k})|)^{\mu}}\dif x \lesssim \left\vert\int_{\Omega}\langle f'_{k}(\sg(u_{k})),\rho\odot\tau_{s,h}^{+}u_{k}\rangle\dif x\right\vert + \frac{1}{k}\|\tau_{s,h}^{-}(\rho^{2}\tau_{s,h}^{+}u_{k})\|_{(\sobo_{0}^{1,\infty})^{*}}, 
\end{align}
with the constants implicit in '$\lesssim$' being uniformly bounded in $k$. Let us briefly explain how to handle the two terms on the right side: As to the first term, we note that $f'_{k}(\sg(u_{k}))$ converges in a suitable sense to the solution of the dual problem associated with \eqref{eq:varprin} (in the sense of convex duality, see \cite{ET}). Since the dual solution $\sigma\in\lebe^{\infty}(\Omega;\R^{n\times n})$ itself belongs to $\sobo_{\locc}^{1,2}(\Omega;\R^{n\times n})$, this regularity can be shown to inherit to $f'_{k}(\sg(u_{k}))$ uniformly in $k$. If we wish to fruitfully use this estimate, we need to suitably bound $\|\rho\tau_{s,h}^{+}u_{k}\|_{\lebe^{2}}$ uniformly in $k$ too. For general $n$, the fractional Sobolev--type embeddings $\bd_{\locc}\hookrightarrow \sobo_{\locc}^{s,n/(n-1+s)}$, $0<s<1$, with the fractional Sobolev spaces $\sobo^{\theta,r}$, $0<\theta<1$ and $r\geq 1$, are optimal, and $n/(n-1+s)=2$ is achieved if and only if $n=2$ and $s=0$. In this case, however, we loose all smoothness information and hence may invoke a condition that is slightly weaker than a local boundedness assumption in the spirit of \eqref{eq:locbound}; namely, we require the so--called \emph{local $\bmo$--assumption}, meaning that for each relatively compact $K\subset\Omega$, the seminorms $\|u_{k}\|_{\bmo(K;\R^{n})}$ are bounded independently of $k$. As shall be demonstrated in \cite{GK} by means of so--called Dorronsoro--type estimates which have been fruitfully used in \cite{KM} by Kristensen \& Mingione in a different context, we have 
\begin{align}\label{eq:embed}
\bd_{\locc}\cap\bmo_{\locc}\hookrightarrow \sobo_{\locc}^{\frac{1}{p}-\varepsilon,p}
\end{align}
for all $1\leq p <\infty$ and suitably small $\varepsilon>0$; the limiting case $\varepsilon=0$ is not even reached in general even if $\bd_{\locc}$ is replaced by the considerbaly smaller space $\sobo_{\locc}^{1,1}$, a fact which has been pointed out by Bourgain, Brezis \& Mironescu in \cite{BBM}, Remark 3. The upshot of this in comparison with the aforementioned embedding without the $\bmo$--side constraint is that although the $\bmo$--condition is not reflected by the first derivatives, it improves both fractional differentiability and the corresponding integrability at a uniform rate. Putting $p=2$ in \eqref{eq:embed}, it is possible to estimate the first term on the right side of \eqref{eq:coercive} by $Ch^{3/2-\varepsilon}$ for any suitably small $\varepsilon>0$. The second term on the right side of \eqref{eq:coercive} can be estimated in the same way, using $\|\Delta_{s,h}^{+}v\|_{(\sobo_{0}^{1,\infty})^{*}}\lesssim \|v\|_{\lebe^{1}}$ together with $\bd_{\locc}\hookrightarrow\sobo_{\locc}^{\theta,1}$ for any $0<\theta<1$. Going back to \eqref{eq:coercive}, we then obtain the uniform bound
\begin{align*}
\int_{\Omega}\left\vert \frac{\rho\tau_{s,h}^{+}\sg(u_{k})}{h^{\frac{3}{4}-\frac{\varepsilon}{2}}}\right\vert^{2}\omega_{k}\dif x := \int_{\Omega}\left\vert \frac{\rho\tau_{s,h}^{+}u_{k}}{h^{\frac{3}{4}-\frac{\varepsilon}{2}}}\right\vert^{2}\frac{1}{(1+|\sg(u_{k})|)^{\mu}}\dif x \leq C
\end{align*}
This is a weighted Nikolski\u{\i}--type estimate for $\sg(u_{k})$. At this stage of the proof, it possible to deduce that the weights $\omega_{k}$ uniformly belong to certain Muckenhoupt classes $A_{p}$ for some $1<p<\infty$ and hence, using embedding results for weighted Nikolski\u{\i} spaces and the theory of singular integrals on Muckenhoupt weighted Lebesgue spaces, it is possible to deduce for a non--empty range of $\mu\in (1,2)$ that the symmetric gradients $\sg(u_{k})$ are uniformly bounded in certain non--weighted $\lebe_{\locc}^{p}$--spaces with $p>1$. By arbitariness of the intially chosen minimising sequence, this implies $u\in\sobo_{\locc}^{1,p}(\Omega;\R^{n})$ by Korn's Inequality and thus establishes the higher integrability of generalised minima of \eqref{eq:varprin} subject to the above local $\bmo$--assumption. In summary, the strategy leads to the following theorem which shall be established in \cite{GK}:
\begin{theorem}\label{thm:main2}
Let $f\in\hold^{2}(\R_{\sym}^{n\times n})$ be a convex integrand of linear growth and $\Omega$ an open and bounded Lipschitz subset of $\R^{n}$. For every $n\in\mathbb{N}$ with $n\geq 2$ there exists a number $1<\mu(n)<2$ and an exponent $p\geq 1$ such that the following holds: If $f$ is $\mu$--elliptic with $\mu\leq\mu(n)$, then any generalised minimiser $u\in\bd(\Omega)\cap\bmo_{\locc}(\Omega;\R^{n})$ of the functional $\mathfrak{F}$ given by \eqref{eq:varprin} belongs to $\sobo_{\locc}^{1,p}(\Omega;\R^{n})$ for some $p\geq 1$. 
\end{theorem}
Let us remark that it is plausible for the preceding result to hold true for a fairly larger range of $\mu$ than described in Theorem \ref{thm:main2}, however, this seems hard to be achieved by use of the above method. In particular, a Sobolev regularity result regarding the limiting case $\mu=3$ (an instance of which is the area integrand $f(\mathbf{Z}):=\sqrt{1+|\mathbf{Z}|^{2}}$) would be desirable. Finally, the strategy outlined above allows to weaken even the local BMO--assumption in view of uniform local $\lebe^{p}$--bounds on the single members of minimising sequences, and shall be addressed in a future publication. 
\section{$\mathbb{A}$--Convex Functionals}\label{sec:Aconvex}
It is a natural extension of the problems and results outlined in the previous sections to replace the symmetric gradient operator by an elliptic differential operator of the form \eqref{eq:differentialoperator}. In turn, within the framework of section \ref{sec:Ornstein}, we aim for existence and regularity properties for minima functionals of the form 
\begin{align*}
\mathfrak{F}[v]:=\int_{\Omega}f(\mathbb{A}[D]u)\dif x
\end{align*}
over Dirichlet classes $\widetilde{\mathscr{D}}:=u_{0}+\sobo_{0}^{\mathbb{A},1}(\Omega)$, where $f\colon W\to\R$ is of linear growth, thus verifying \eqref{eq:lingrowth} with the obvious modifications. The spaces $\sobo^{\mathbb{A},1}(\Omega)$ are defined similarly to $\sobo^{1,1}$ or $\ld$, namely, we say that a measurable map $v\colon\Omega\to V$ belongs to $\sobo^{\mathbb{A},1}(\Omega)$ if and only if $v\in\lebe^{1}(\Omega;V)$ and the weak differential expression $\mathbb{A}[D]u$ belongs to $\lebe^{1}(\Omega;W)$. Equipped with the canonical norm, one defines $\sobo_{0}^{\mathbb{A},1}(\Omega)$ as the closure of $\hold_{c}^{1}(\Omega;V)$ with respect to the norm topology. Similarly, we define the space of functions of bounded $\mathbb{A}[D]$--variation $\bv_{\mathbb{A}}(\Omega)$ to consist of all $v\in\lebe^{1}(\Omega;V)$ for which the distributional differential expression $\mathbb{A}[D]v$ can be represented by a finite $W$--valued Radon measure. Based on the linear growth assumption, the direct method in conjunction with the obvious changes of section \ref{sec:relaxed} would lead to a satisfactory existence theory within the class of functions bounded $\mathbb{A}[D]$--variation \emph{provided} the trace operator for such function spaces would be well--understood. This, however, is not the case by now: For instance, even if $\mathbb{A}[D]$ is an elliptic operator in the sense of section \ref{sec:Ornstein}, elements of $\bv_{\mathbb{A}}$ or $\sobo^{\mathbb{A},1}$ do not neccesarily possess traces in $\lebe^{1}(\partial\Omega;V)$ for arbitrary Lipschitz domains $\Omega\subset\R^{n}$: 
\begin{remark}
For $n=2$, let the trace--free symmetric gradient operator defined by $\sg^{D}(v):=\sg(v)-\frac{1}{2}\di(v)\mathbbm{1}_{2\times 2}$ with the $(2\times 2)$--unit matrix. In this situation, $\sg^{D}$ is an elliptic, linear, homogeneous first order differential operator  on $\R^{2}$ from $V=\R^{2}$ to $W=\R^{2\times 2}$. Identifying $\R^{2}\cong\mathbb{C}$, it is easy to see that $\ker(\sg^{D})$ contains all holomorphic functions. To argue that elements of $\sobo^{\sg^{D},1}(\ball(0,1);\R^{2})(\simeq \sobo^{\sg^{D},1}(\mathbb{D};\mathbb{C}))$ with the unit disk $\mathbb{D}\subset\mathbb{C}$ do not have traces in $\lebe^{1}(\partial\ball(0,1);\R^{2})(\simeq \lebe^{1}(\partial\mathbb{D};\mathbb{C}))$, consider $f\colon\mathbb{D}\ni z \mapsto 1/(z-1)\in\mathbb{C}$. Then $f\in\lebe^{1}(\mathbb{D};\mathbb{C})$, is holomorphic and thus belongs to $\sobo^{\sg^{D},1}(\mathbb{D};\mathbb{C})$ whereas it is easy to see that $\int_{\partial\mathbb{D}}|f(z)|\dif z = \infty$. 
\end{remark}
The previous example is due to Fuchs \& Repin \cite{FR}, and motivates the characterisation of all $\mathbb{A}[D]$ such that the corresponding spaces $\sobo^{\mathbb{A},1}(\mathbb{D};\mathbb{C})$ possess trace space $\lebe^{1}(\partial\Omega;V)$ at least for the large class of bounded Lipschitz subsets $\Omega$ of $\R^{n}$. For such operators, the Sobolev regularity result, Theorem \ref{thm:main2} is easily shown to hold true as its proof does not use the specific structure of the symmetric gradient operator. In this respect, it is important to note that the techniques available in the literature -- so for instance Strang \& Temam's approach \cite{ST} in the case of $\bd$ -- which are taylored for the symmetric gradient case, do not apply to the setting of arbitrary elliptic differential operators of the form \eqref{eq:differentialoperator} without substantial modifications. We hope to succesfully tackle this problem in a future publication.

\end{document}